

\font\sectionfont=cmbx10 scaled\magstep1
\def\titlea#1{\vskip0pt plus.3\vsize\penalty-75
    \vskip0pt plus -.3\vsize\bigskip\bigskip
    \noindent{\sectionfont #1}\nobreak\smallskip\noindent}

\def\claim#1#2{\vskip.1in\medbreak\noindent{\bf #1.} {\sl #2}\par
    \ifdim\lastskip<\medskipamount\removelastskip\penalty55\medskip\fi}
\def\rmclaim#1#2{\vskip.1in\medbreak\noindent{\bf #1.} {#2}\par
    \ifdim\lastskip<\medskipamount\removelastskip\penalty55\medskip\fi}
\def\beglemma#1#2\endlemma{\claim{#1 Lemma}{#2}}
\def\begdefinition#1#2\enddefinition{\claim{#1 Definition}{#2}}
\def\begtheorem#1#2\endtheorem{\claim{#1 Theorem}{#2}}
\def\begcorollary#1#2\endcorollary{\claim{#1 Corollary}{#2}}
\def\begremark#1#2\endremark{\rmclaim{#1 Remark}{#2}}
\def\begproposition#1#2\endproposition{\claim{#1 Proposition}{#2}}
\def\begassumption#1#2\endassumption{\claim{#1}{#2}}
\def\begProof{\noindent{\bf Proof.}\quad}
\def\begProofof#1{\medskip\noindent{\bf Proof of #1.}\quad}
\def\square{\hbox{$\sqcap\!\!\!\!\sqcup$}}
\def\endProof{\hfill\square\par
    \ifdim\lastskip<\medskipamount\removelastskip\penalty55\medskip\fi}
\newcount\FNOTcount \FNOTcount=1
\def\numfnot{\number\FNOTcount}
\def\addfnot{\global\advance\FNOTcount by 1}
\def\fonote#1{\footnote{$^\numfnot$}{#1}\addfnot}
\def\acknow#1{{\bf Acknowledgements.}\quad #1}

\def\Co{{\bf C}} 
\def\Re{{\bf R}} 
\def\To{{\bf T}} 
\def\Ze{{\bf Z}} 

\def\A{{\cal A}}
\def\a{\alpha}
\def\B{{\cal B}}

\def\C{{\cal C}}

\def\D{\Delta}
\def\cD{{\cal D}}
\def\f{\varphi}
\def\F{{\cal F}}

\def\G{\Gamma}
\def\bG{{\bf G}}
\def\H{{\cal H}}
\def\k{\kappa}
\def\K{{\cal K}}
\def\bK{{\bf K}}
\def\l{\lambda}
\def\L{\Lambda}
\def\m{\mu}
\def\n{\nu}

\def\O{{\cal O}}
\def\p{\pi}
\def\Q{\Omega}
\def\r{\rho}
\def\s{\sigma}
\def\S{{\cal S}}

\def\T{\Theta}
\def\th{\vartheta}
\def\W{{\cal W}}

\def\Z{{\cal Z}}
\def\Lpo{{\cal L}_+^\uparrow}           
\def\Pp{{\cal P}_+}                     
\def\Ppo{{\cal P}_+^\uparrow}           
\def\Sp{\widetilde{\cal P}_+}           
\def\Spo{\widetilde{\cal P}_+^\uparrow} 
\def\imply{\Rightarrow}
\def\Ind#1#2{{\rm Ind}_{#1}^{#2}}

\def\npi{\par}


\def\refno#1#2{\item{[#1]}{#2}}
\def\begref#1#2{\titlea{#1}}
\def\endref{}

\newcount\REFcount \REFcount=1
\def\numref{\number\REFcount}
\def\addref{\global\advance\REFcount by 1}
\def\wdef#1#2{\expandafter\xdef\csname#1\endcsname{#2}}
\def\wdch#1#2#3{\ifundef{#1#2}\wdef{#1#2}{#3}
    \else\write16{!!doubly defined#1,#2}\fi}
\def\wval#1{\csname#1\endcsname}
\def\ifundef#1{\expandafter\ifx\csname#1\endcsname\relax}

\def\autonumref{
    \def\rfr(##1){\wdef{q##1}{yes}\ifundef{r##1}$\diamondsuit$##1
        \write16{!!ref ##1 was never defined!!}\else\wval{r##1}\fi}
    \def\REF(##1)##2\endREF{\wdch{r}{##1}{\numref}\addref}\REFERENCES
    \def\references{
        \def\REF(####1)####2\endREF{
            \ifundef{q####1}\write16{!!ref. [####1] was never quoted!!}\fi
            \refno{\rfr(####1)}####2}
        \begref{References}{99}\REFERENCES\endref}}

\def\createbibfile{
    \def\REF(##1)##2\endREF{\wdch{r}{##1}{defined}}\REFERENCES
    \def\rfr(##1){\wdef{q##1}{quoted}\ifundef{r##1}$\diamondsuit$##1
    \write15{\string \REF(##1) [##1] to be inserted \csname endREF\endcsname}
        \write16{!!ref ##1 was never defined!!}\else##1\fi}
    \def\references{
        \def\REF(####1)####2\endREF{
            \ifundef{q####1}
            \else\refno{####1}{####2}
    \write15{\string \REF(####1)####2\csname endREF\endcsname}\fi}
    \begref{References}{99}\REFERENCES\endref}}

 \def\REFERENCES{
 \REF(AHKT1)Araki, H., Haag, R., Kastler, D., Takesaki, M.: ``{\it
Extensions of KMS states and chemical potential.}" Commun.
Math. Phys. {\bf 53}, 97-134 (1977)
 \endREF
 \REF(BiWi1)Bisognano, J., Wichmann, E.: ``{\it On the duality
condition for a Hermitian scalar field.}" J. Math. Phys. {\bf
16}, 985-1007 (1975)
 \endREF
 \REF(BiWi2)Bisognano, J., Wichmann, E.: ``{\it On the duality
condition for quantum fields.}" J. Math. Phys. {\bf
17}, 303-321 (1976)
 \endREF
 \REF(BjDr1)Bjorken, J.D., Drell, S.D.: {\it Relativistic
Quantum Fields.} New York: McGraw-Hill, 1965
 \endREF
 \REF(Borc2)Borchers, H.J.: ``{\it Local rings and the
connection between spin and statistics.}" Commun. Math. Phys. {\bf 1},
281-307 (1965)
 \endREF
 \REF(Borc1)Borchers, H.J.: ``{\it The CPT theorem in
two-dimensional theories of local observables.}" Commun. Math.
Phys. {\bf 143}, 315 (1992)
 \endREF
 \REF(Borc3)Borchers, H.J.: ``{\it On the converse of the Reeh-Schlieder
theorem.}" Commun. Math. Phys. {\bf 10}, 269-273 (1968)
 \endREF
 \REF(BGL1)Brunetti, R., Guido, D., Longo, R.: ``{\it Modular
structure and duality in conformal quantum field theory.}"
Commun. Math. Phys. {\bf 156}, 201-219 (1993)
 \endREF
 \REF(BGL2)Brunetti, R., Guido, D., Longo, R.: ``{\it Group
cohomology, modular theory and space-time symmetries.}"
Rev. Math. Phys. {\bf 7}, 57-71 (1994)
 \endREF
 \REF(BuEp1)Buchholz, D., Epstein, H.: ``{\it Spin and statistics
of quantum topological charges.}" Fizika {\bf 3},
329-343 (1985)
 \endREF
 \REF(BuFr1)Buchholz, D., Fredenhagen, K.: ``{\it Locality and
structure of particle states.}" Commun. Math. Phys. {\bf 84},
 1-54 (1982)
 \endREF
 \REF(BuSu1)Buchholz, D., Summers, S.J.: ``{\it An algebraic
characterization of vacuum states in Minkowski space.}" Commun. Math. Phys.
{\bf 155}, 442-458 (1993)
 \endREF
 \REF(Burgoyne)Burgoyne, N.: ``{\it On the connection of spin with
statistics.}" Nuovo Cimento {\bf 8}, 807 (1958)
 \endREF
 \REF(Dell'Antonio)Dell'Antonio, G.F.: ``{\it On the connection of spin with
statistics.}" Ann. Phys. {\bf 16}, 153 (1961)
 \endREF
 \REF(DHR1)Doplicher, S., Haag, R., Roberts, J.E.: ``{\it Local
observables and particle statistics I.}" Commun. Math. Phys.
{\bf 23}, 199-230 (1971)
 \endREF
 \REF(DHR2)Doplicher, S., Haag, R., Roberts, J.E.: ``{\it Local
observables and particle statistics II.}" Commun. Math. Phys.
{\bf 35}, 49-85 (1974)
 \endREF
 \REF(DoLo1)Doplicher, S., Longo, R.: ``{\it Standard and split
inclusions of von Neumann algebras.}" 	Invent. Math. {\bf 73}, 493-536
(1984)
 \endREF
 \REF(DoRo)Doplicher, S., Roberts, J.E.: ``{\it Why there is a
 field algebra with a compact gauge group describing the
 superselection structure in particle physics.}" Commun. Math.
 Phys. {\bf 131}, 51-107 (1990)
 \endREF
 \REF(Dyson)Dyson, F.J.: ``{\it On the connection of weak local
commutativity and regularity of Wightman functions.}" Phys. Rev.
{\bf 110}, 579 (1958)
 \endREF
 \REF(Fie)Fierz, M.: ``{\it \"Uber die relativische Theorie
kr\"aftfreier Teilchen mit beliebigem spin.}" Helv. Phys. Acta
{\bf 12}, 3 (1939)
 \endREF
 \REF(Epst1)Epstein, H.: ``{\it CTP invariance in a theory of
local observables.}" J. Math. Phys. {\bf 8}, 750 (1967)
 \endREF
 \REF(GL1)Guido, D., Longo, R.: ``{\it Relativistic invariance
and charge conjugation in quantum field theory.}" Commun. Math.
Phys. {\bf 148}, 521-551 (1992)
 \endREF
 \REF(GL3)Guido, D., Longo, R.: ``{\it The Conformal Spin and Statistics
theorem.}". Preprint.
 \npi
Brunetti, R., Guido, D., Longo, R.: in preparation.
 \endREF
  \REF(Haag1)Haag, R.: {\it Local Quantum Physics.} Berlin Heidelberg:
Springer Verlag, 1992
 \endREF
 \REF(HHW)Haag, R., Hugenoltz, N.M., Winnink, M.: ``{\it On the equilibrium
states in quantum statistical mechanics.}" Commun. Math.
Phys. {\bf 5}, 215 (1967)
 \endREF
 \REF(Hawk1) Hawking, S.W.: ``{\it Particle creation by black holes.}"
Commun. Math. Phys. {\bf 43}, 199 (1975)
 \endREF
 \REF(HaKa1)Haag, R., Kastler, D.: ``{\it An
algebraic approach  to Quantum Field Theory.}" J. Math. Phys. {\bf 5},
848-861 (1964)
 \endREF
 \REF(Jost1)Jost, R.: {\it The general theory of Quantized Fields.}
 Providence RI: Am. Math. Soc., 1965
 \endREF
 \REF(Jos57)Jost, R.: ``{\it Eine Bemerkung zu CTP Theorem.}"
Helv. Phys. Acta {\bf 30}, 409 (1957)
 \endREF
 \REF(Kay) Kay B.S., ``{\it A Uniqueness Result for Quasi-Free KMS States}''
Helv. Phys. Acta {\bf 58}, 1017-1029 (1985)
 \endREF
 \REF(Kuck1)Kuckert, B.: ``{\it  PCT und Kovarianz als modulare Strukturen.}"
Diploma Thesis, Hamburg, in preparation
 \endREF
 \REF(Lipsman)Lipsman, R.L.: ``{\it Group representations.}" Lect. Notes
in Math. {\bf 388}, New York--Heidelberg--Berlin: Springer Verlag, 1974
 \endREF
 \REF(Land1)Landau, L.: ``{\it On local functions of fields.}" Commun. Math.
Phys. {\it 39}, 49-62 (1974)
\endREF
 \REF(Long1)Longo, R.: ``{\it Index of subfactors and statistics
of quantum fields.}" $I$ Commun. Math. Phys. {\bf 126}, 217-247 (1989);
$II$ {\bf 130}, 285-309 (1990)
 \endREF
 \REF(Lud54)L\"uders, G.: ``{\it Vertaugschungsrelationen zwischen
verschiedenen Feldern.}" Z. Na\-tur\-forsch. {\bf 13a}, 254 (1958)
 \endREF
 \REF(LudersAndZumino)L\"uders, G., Zumino, B.: ``{\it Connection between spin
and  statistics.}" Phys. Rev. {\bf 110}, 1450 (1958)
 \endREF
 \REF(Pau)Pauli, W.: ``{\it On the connection between spin and statistics.}"
Phys. Rev. {\bf 58}, 716 (1940)
 \endREF
 \REF(Pau55)Pauli, W.: ``{\it Exclusion principle, Lorentz group and
reflection of space time and charge.}" In: Niels Bohr and the Development of
Physics, W. Pauli (ed.) New York: Pergamon Press, 1955
 \endREF
 \REF(Schwinger)Schwinger, J.: ``{\it On the theory of quantized fields I.}"
Phys. Rev. {\bf 82}, 914 (1951)
 \endREF
 \REF(Sewell)Sewell, G.L.: ``{\it Relativity of temperature and Hawking
effect.}" Phys. Lett. {\bf79A}, 23 (1980)
 \endREF
 \REF(StZs1)Str\u atil\u a, S., Zsido, L.: {\it Lectures on von~Neumann
algebras.} England: Abacus Press, 1979
 \endREF
 \REF(Streater)Streater, R.F.: {\it Local fields with the wrong
connection between Spin and Statistics.} Commun. Math. Phys.
{\bf 5}, 88-96 (1967)
 \endREF
 \REF(StWi1)Streater, R.F., Wightman, A.S.: {\it PCT, Spin and Statistics,
and all that.} Reading (MA): Benjamin, 1964
 \endREF
 \REF(Take1)Takesaki, M.: ``{\it Tomita theory of modular Hilbert
algebras.}" Lect. Notes in Math. {\bf 128}, New~York--Heidelberg--Berlin:
Springer Verlag, 1970
 \endREF
 \REF(Unru1) Unruh, W.G.: ``{\it Notes on black hole evaporation.}"
Phys. Rev. {\bf D14}, 870 (1976)
 \endREF
 \REF(Wigh)Wightman, A.S.: ``{\it Quantum field theory in terms
of vacuum expectation values.}" Phys. Rev. {\bf 101}, 860 (1956)
 \endREF
 \REF(Wies1)Wiesbrock, H.V.: ``{\it A comment on a recent work of
Borchers.}" Lett. Math. Phys. {\bf 25}, 157-159 (1992)
 \endREF
\REF(Yngv1)Yngvason, J.: ``{\it A note on essential
duality.}" Lett. Math. Phys. {\bf31}, 127-141 (1994)
\endREF}

\autonumref
%

\topskip3.cm
\font\ftitle=cmbx12 scaled\magstep1
\vskip2truecm
\centerline{\ftitle An Algebraic Spin and Statistics Theorem}

 \bigskip\bigskip

\centerline{Daniele Guido$^1$\footnote{$^*$}
{ Supported in part by MURST and
CNR-GNAFA.} and Roberto Longo$^{1,2*}$}
\footnote{}{E-mail:\ guido@mat.utovrm.it,
longo@mat.utovrm.it }
 \vskip1.truecm \item{$(^1)$}
 Dipartimento di Matematica, Universit\`a di
Roma ``Tor Vergata'' \par via
della Ricerca Scientifica,
I--00133 Roma, Italia.
\item{$(^2)$} Centro Linceo Interdisciplinare,
Accademia Nazionale dei Lincei
\par via della Lungara 10, I--00165
Roma, Italia
\bigskip
{\it Dedicated to Hans Borchers on the occasion of his seventieth
birthday}
\bigskip
\bigskip
\hfill To appear in Commun. Math. Phys. (1995).
\bigskip
\bigskip
{\bf Abstract:} A spin-statistics theorem and a PCT theorem are obtained
in the context of the superselection sectors in  Quantum Field Theory
on a 4-dimensional space-time.
 Our main assumption is the requirement that the modular groups of the
von Neumann algebras of local observables associated with wedge regions act
geometrically as pure Lorentz transformations.
 Such a property, satisfied by the local algebras generated
by Wightman fields because of the Bisognano-Wichmann theorem, is regarded as
a natural primitive assumption.

\vfill\eject
\topskip0.cm

\titlea{Introduction}

In this paper we shall reconsider from an intrinsic point of view two well
known fascinating theorems in Quantum Field Theory: the PCT theorem
and the Spin and Statistics theorem.

Both of them have a long history, see [\rfr(Jost1),\rfr(StWi1)]. The spin and
statistics theorem first appeared in the context of free fields in the work of
Fierz [\rfr(Fie)] and Pauli [\rfr(Pau)]: one cannot second quantize particles
with integer spin by anticommuting fields, i.e. fields obeying Fermi
statistics,
nor particles with half-integer spin by local fields, i.e. fields obeying Bose
statistics.

The PCT theorem originated in [\rfr(Lud54)] as a relation between the
existence of the space-inversion symmetry P and the existence of the product of
the charge and the time-inversion symmetry CT. Pauli proved in
[\rfr(Pau55)] that PCT is always a symmetry of Lorentz invariant field
equations.

It was a success of the Wightman axiomatic approach [\rfr(Wigh)] to establish
model independent results: the connection between spin and statistics was
obtained by Burgoyne [\rfr(Burgoyne)], see also
[\rfr(LudersAndZumino),\rfr(Dell'Antonio)], and a PCT theorem by Jost
[\rfr(Jos57)], see also [\rfr(Dyson), \rfr(Schwinger)], both relying on the
general holomorphic properties of the $n$-point functions. A  spin and
statistics theorem in the algebraic approach [\rfr(HaKa1)], see also
[\rfr(Borc2)], was later given by Epstein [\rfr(Epst1)] and has a version for
(Doplicher-Haag-Roberts) DHR superselection sectors [\rfr(DHR2)] and for
more general topological charges [\rfr(BuEp1)].

All these approaches heavily rely on arguments of analytic continuation, whose
nature give some mysterious effectiveness to the results. Moreover they make
use of certain detailed structures, either because they deal with Wightman
tempered distributions or because  they treat the case of finite mass
degeneracy, where a supeselection sector has to contain only finitely many
particles of the same type and all of them are assumed to have strictly
positive mass.

The approach to Quantum Field Theory by local observable algebras
[\rfr(HaKa1)] suggests however that a PCT symmetry and spin-statistics
correspondence should be intrinsically associated with the net of local
algebras and manifest itself as the consequence of the locality
principle\fonote{The spin-statistics relation depends on  sharp locality. In
the second quantization of Bose particles by anti-commuting fields, microscopic
causality is still asymptotically present and its violation is sizeable only
at distances comparable to the Compton wave length [\rfr(BjDr1)]. This is
perhaps an indication that the spin-statistics relation might be different in
contexts like Quantum Gravity where a sharp causality principle does not
occur.}.

{}From the mathematical point of view the spin-statistics correspondence is a
relation between two  quantities of different nature, the univalence and
the statistics phase, and one is led to tie up these concepts on general
grounds, somehow in the spirit of an index theorem.

We shall establish both a PCT and a spin-statistics theorem in the following
general context.

Let $\O\to\A(\O)$ be a net of von~Neumann algebras on a Hilbert space
$\H$, i.e. an inclusion preserving association  between regions
$\O$ in the four-dimensional Minkowski space and von Neumann
algebras of local observables, that we assume here to be irreducible.

 We make the  following assumptions.

 {\it Locality.} If $\O_1$ and $\O_2$ are space-like separated regions, then
$\A(\O_1)$ and $\A(\O_2)$ commute elementwise.

{\it Modular covariance.} There is a vector $\Q\in\H$, the vacuum
vector,  cyclic for the algebras $\A(W)$ associated with all wedge region $W$
in
the Minkowski space, such that
 $$
 \D_W^{it}\A(\O)\D_W^{-it}=\A(\L_W(t)\O)\ ,\qquad t\in\Re
 $$
 where $\O$ is any region, $\D$ is the Tomita-Takesaki modular operator
[\rfr(Take1),\rfr(StZs1)] associated with $(\A(W),\Q)$ and $\L_W$ is the
one-parameter rescaled group of pure Lorentz transformations preserving $W$.

{\it Reeh-Schlieder property.} The vacuum vector $\Q\in\H$ is
also cyclic for the
algebras $\A(\S)$ associated with all space-like cones $\S$.

Locality is the well-known expression of Einstein causality and we do not
dwell on it. The Reeh Schlieder property is
known to hold for the vacuum vector in a Poincar\'e covariant theory as a
consequence of the positivity of the energy and the weak additivity assumption
for the local algebras [\rfr(Borc3)].

Modular covariance  needs however further comments. Postponing for a while
the justification for such an assumption, we recall that this entails  the net
to
be covariant with respect to the universal covering $\Spo$ of the Poincar\'e
group $\Ppo$, with positive energy [\rfr(BGL2)]. Indeed we shall prove here
that
it is actually covariant with respect to $\Ppo$ as a special case of our
general
Spin and Statistics theorem and taking to completion our previous work.
Therefore  modular covariance is a way to intrinsically encode the Poincar\'e
covariance property   in the net structure,  providing a
canonical representation of the Poincar\'e group $\Ppo$ (cf. also
[\rfr(BuSu1)]).

Let now $\r$ be a superselection sector of $\A$ in the sense of
Doplicher-Haag-Roberts [\rfr(DHR1)] or more generally of Buchholz-Fredenhagen
[\rfr(BuFr1)].

An index-statistics relation [\rfr(Long1)] shows that
 $$
 {\rm Ind}(\r)=d(\r)^2
 $$
 where ${\rm Ind}(\r)$ is the Jones index of $\r$ and $d(\r)$ is the DHR
statistical dimension,  namely
 $$
 d(\r)=|\l_\r|^{-1}
 $$
 where the statistical parameter $\l_\r\in\Re$ classifies the statistics in
$3+1$ space-time dimensions [\rfr(DHR1)].

Therefore the index is an intrinsic quantity that determines the statistics up
to the Fermi-Bose alternative, i.e. the sign of $\l_\r$.

On the other hand, the Poincar\'e representation in the vacuum sector being
fixed by modular covariance, the representation of $\Spo$ associated with a
covariant irreducible sector $\r$ is uniquely determined, therefore the
univalence (integer or half-integer spin alternative) is intrinsically
associated with $\r$. Since  $\r$ is automatically $\Spo$-covariant if
$d(\r)<\infty$ (assuming a regularity property for the net [\rfr(GL1)]), it is
natural to expect a general algebraic Spin and Statistics theorem connecting
these two intrinsic quantities for any sector with finite statistics.

Our result in this respect will in fact show that on these general grounds
 $$
 {\rm sign}(\l_\r)=U_\r(2\pi)
 $$
 where $U_\r$ is the representation of $\Spo$ in the sector $\r$ and
$U_\r(2\pi)$ denotes the corresponding rotation by $2\pi$.

 Modular covariance also implies  that the anti-unitary involution $\T$,
definable by the modular theory according to the Bisognano-Wichmann
prescription [\rfr(BiWi2)], implements a complete space-time reflection. As
shown in [\rfr(GL1)], this entails that $\T$ intertwines a sector with its
conjugate. We therefore  obtain a PCT symmetry.

We come  now  back to the origin of the modular covariance
property. Its main justification certainly comes from the
Bisognano-Wichmann theorem [\rfr(BiWi1),\rfr(BiWi2)], to the effect that this
property holds if the local algebras are constructed from Wightman fields.

An algebraic version of the Bisognano-Wichmann theorem
does not exist yet, except in the case of conformal theories where it holds in
full generality [\rfr(BGL1)].
However a theorem of Borchers [\rfr(Borc1)] shows  part of the geometric
properties of the modular group for wedge regions to be always present and in
particular every $1+1$ dimensional Poincar\'e covariant net satisfying
essential
duality has the modular covariance property.

At the present time no counter-example to modular covariance is known to
exist within Poincar\'e covariant theories (see however
[\rfr(Land1),\rfr(Yngv1)]). There are nevertheless counter-examples to the
spin-statistics theorem [\rfr(Streater)]: these are constructed by infinite
multiplicity fields where the Poincar\'e group representation is not unique. It
turns out that the wrong connection between spin and statistics depends on the
wrong choice of the Poincar\'e group representation, while our canonical choice
for the latter has the desired properties. We remark that an intrinsic way to
eliminate pathological examples of the above kind comes by requiring the split
property [\rfr(DoLo1)]; this indeed implies the uniqueness of the Poincar\'e
group representation [\rfr(BGL1)] and we propose it as a natural candidate for
a
derivation of the modular covariance property by first principles.

On the physical side modular covariance manifests an interesting analogy with
the
Unruh effect [\rfr(Unru1)] and with the Hawking black hole thermal radiation
[\rfr(Hawk1)], as first noticed by Sewell [\rfr(Sewell)]. We sketch
the essential ideas, see also [\rfr(Haag1)]. As is known the modular group of a
von Neumann algebra with respect to a given state is characterized by the
Kubo-Martin-Schwinger condition [\rfr(Take1)] and, on the other hand, the KMS
condition is peculiar of thermal equilibrium states in Statistical Mechanics
[\rfr(HHW)]. By the Bisognano-Wichmann theorem the boosts satify the KMS
condition with respect to the vacuum, as automorphisms of the   von Neumann
algebra of the corresponding wedge $W$ and, on the other hand, the orbits of
the
boosts are the trajectories of a uniformly accelerated motion for which the
"Rindler universe" $W$ is a natural horizon; the equivalence principle in
Relativity Theory then allows an interpretation of the thermal outcome as a
gravitational effect. On this basis Haag has proposed long ago to derive the
Bisognano-Wichmann theorem.

The role of the modular covariance assumption may be also understood by its
consequences. Among other things, it implies the positivity of the energy for
the constructed Poincar\'e group representation [\rfr(Wies1),\rfr(BGL2)]. As is
known the positivity of the energy is lost on a curved space-time, and the
modular covariance seems to be the appropriate substitute in this case.
 Moreover, as already mentioned, it gives rise to the KMS
condition, namely an analytic continuation property. It turns out that
this analytic aspect of the modular covariance assumption incorporates all the
 holomorphic properties present in Quantum Field Theory. But, as a matter of
fact, the modular group is an algebraic object, a manifestation of the
$^*$-operation, thus providing us with an algebraic approach to our problems.

We pass now to a description of the methods of our work. This paper relies on
the modular theory of Tomita and Takesaki and on an analysis by the unitary
representations of ${\rm SL}(2,\Re)$.

We shall find a key relation arising from the comparison of the modular
groups of different algebras, and we shall regard it as an identity concerning
operators in the space of a representation of ${\rm SL}(2,\Re)$, because of the
well known fact that the 2+1 dimensional Lorentz group is isomorphic to
${\rm SL}(2,\Re)/\{1,-1\}$. Section~1 contains the proof of this identity by
the Mackey machine of the induced representations (e.g. [\rfr(Lipsman)]) and a
free field verification.

In Sect.~2 the PCT theorem and the Spin and Statistics relation are proven
in the context of the field algebras,\fonote{We have  recently been informed by
Kuckert of an  independent complementary  analysis based on  assumptions of
weak
geometric type for the modular conjugation [\rfr(Kuck1)].} where the formalism
is close to  the classical formulation.
 Then, in Sect.~3, we obtain our result  in the context  of local
observables. This is done by rephrasing the statements in terms of the
Doplicher-Roberts field algebra [\rfr(DoRo)]. This last step has certain
pedagogical advantages, but has to be avoided in order
to extend our work to more
general settings where  the field algebra does not exist.

In forthcoming papers [\rfr(GL3)] we shall indeed provide a more intrinsic
approach in terms of local observables only, that will cover low dimensional
and
conformal theories  in particular. The general picture will be clarified by
examples.

\titlea{1. An Identity for Operators Associated with
 Representations of ${\rm SL}(2,\Re)$}

 Let us consider two one-parameter subgroups of ${\rm SL}(2,\Re)$
 $$
\mu(t)\equiv\left (\matrix{\cosh\pi t&-\sinh\pi t\cr-\sinh\pi t&\cosh\pi
t\cr}\right)\ ,\quad
\nu(t)\equiv\left (\matrix{e^{-\pi t}&0\cr0&e^{\pi t}\cr}\right)\ .
 $$

If $U$ is a unitary representation of ${\rm SL}(2,\Re)$ on a Hilbert space
$\H$ we look at the corresponding selfadjoint infinitesimal generators
 $$
H = H_U = i{d\over dt}U(\mu(t))|_{t=0}\ ,
 $$
 $$
K = K_U = i{d\over dt}U(\nu(t))|_{t=0}\ .
 $$

 We shall denote by $\bG$ the group ${\rm PSL}(2,\Re)$ given by the quotient of
${\rm SL}(2,\Re)$ by its center $\{-1,1\}$ and by $\widetilde\bG$ the universal
covering of $\bG$, which is of course the universal covering of ${\rm
SL}(2,\Re)$
too.
 We take the same definition for $H$ and $K$ in the case of a unitary
representation $U$ of the universal covering group $\widetilde\bG$.

 We will consider the following properties for a representation $U$ of
$\widetilde\bG$:
\item{$(i)$} The operator $T_t\equiv e^{{ 1\over 2}K}e^{itH}e^{-{1\over 2}K}$
is densely defined for all $t\in\Re$.
 \item{$(ii)$} $T_t$ is densely defined and
$$
e^{{ 1\over 2}K}e^{itH}e^{-{1\over 2}K}\subset e^{-itH}
 $$
 where the symbol $\subset$ denotes the extension of operators.
 These properties refer to a representation $U$, but we omit the symbol
$U$ when no confusion arises.
 \begtheorem{1.1.} $(a)$ Property~$(ii)$ holds for the regular
representation $\l$ of $\bG$.
 \item{$(b)$} Property $(i)$ implies Property $(ii)$.
 \endtheorem
 In order to prove this theorem, we first observe that it is enough to check
Property~$(ii)$ on dense sets of vectors, not necessarily on the full domain of
$T_t$.

\beglemma{1.2.} Let us assume that, for each real $t$, there is a dense subset
${\cal D}_t$ of the domain of $T_t$ such that $T_t|_{{\cal D}_t}\subset
e^{-itH}$. Then Property~$(ii)$ holds for the given representation.
\endlemma
\begProof Note first that the matrix
 $\left(\matrix{0&1\cr-1&0\cr}\right )\in {\rm SL}(2,\Re)$
 conjugates $\mu(t)$ with $\mu(-t)$ and $\nu(t)$ with $\nu(-t)$
therefore the assumption of the lemma remains true if we replace $K$ with
$-K$ and $H$ with $-H$, in particular $e^{ -{1\over 2}K}e^{itH}e^{{1\over 2}K}$
is densely defined. Now
 $$
T_t^*\supset e^{ -{1\over 2}K}e^{-itH}e^{{1\over 2}K}
 $$
 therefore $T_t^*$ is densely defined and $T_t$ is closable. Since $T_t\xi
= e^{-itH}\xi$ for all $\xi$ in a dense set and $e^{-itH}$ is bounded,
the equality $T_t\xi=e^{-itH}\xi$ must hold for all $\xi$ in the domain of
$T_t$.
 \endProof
 \begcorollary{1.3.}Let $U_1,U_2$ and $U$ be unitary representations of
$\widetilde\bG$.
\medskip
 \item{$(a)$} If Property~$(i)$, resp. $(ii)$, holds for both $U_1$ and $U_2$
then it holds for $U_1\otimes U_2$.
 \item{$(b)$} If Property~$(i)$ holds for both $U_1$ and $U_2$ and
Property~$(ii)$ holds for $U_1\otimes U_2$, then Property~$(ii)$ holds for both
$U_1$ and $U_2$.
 \item{$(c)$}If $U = \int^{\oplus}U(\lambda)dm(\lambda)$ is a direct integral
decomposition of $U$, then Property~$(i)$, resp. $(ii)$, holds for $U$ iff it
holds for $U(\l)$, for $m$-almost all $\lambda$.
\medskip
 \endcorollary
 \begProof Part $(a)$ and $(c)$ are immediate by Lemma~1.2 since one can check
Property~(1.1) on natural dense sets.

 If Property~$(ii)$ holds for $U_1\otimes U_2$ then it holds for $U_1$ and
$U_2$
up to a constant, namely, considering for example the representation $U_1$,
there exists a phase $z(t)$ such that $T_t\subset z(t)e^{-itH}$
and $T_t$ is densely defined. Equivalently $e^{itH}e^{-{1\over2}K}$ is densely
defined and extended by $z(t)e^{-{1\over2}K}e^{-1tH}$ and this implies that
$z(t)$ is a one-dimensional character.

Of course $z(\cdot)$ remains unchanged if we replace $\mu$ and $\nu$ by a pair
of conjugate one-parameter subgroups. As in the proof of Lemma~1.2 we may thus
replace $\mu(t)$ by $\mu(-t)$ and thus $z(t)$ by $z(-t)$, hence $z(t)=z(-t)=1$
and the proof is complete.
 \endProof
 We need now to verify Property~$(ii)$ in some specific representation.
To this end recall that [\rfr(StWi1)], if $W_i$ is the wedge in the
$3$-dimensional space-time along the axis $x_i$, $i=1,2$, and $\L_i(t)$ is the
associated one-parameter group of pure Lorentz transformations
 (see Sect.~2), there is an isomorphism of ${\rm PSL}(2,\Re)$ with the
$2+1$-dimensional Lorentz group $\Lpo(3)$ determined by
 $$
\m(t)\to\L_1(t)\quad{\rm and}\quad
\n(t)\to\L_2(t) .
$$
 Accordingly,  we shall identify $\bG$ with a subgroup of the 2+1-dimensional
Poincar\'e group $\Ppo(3)$.
 \beglemma{1.4.} Let $V\equiv V_{m,0}$ be the positive energy representation of
$\Ppo(3)$ of spin 0 and mass $m>0$. Then Property~$(ii)$ holds for the
restriction $U\equiv V|_{\bG}$ of $V$ to $\bG$.
 \endlemma
 \begProof As is known, $V$ extends to a (anti-)representation of the
proper
Poincar\'e group $\Pp(3)$, namely there exists a anti-unitary involution
$\Theta$, on the same Hilbert space, that commutes with $U$ and
implements the change of sign on the translation operators in any space-time
direction.

 By the one-particle version of the Bisognano Wichmann theorem
(which follows of course from the Bisognano-Wichmann theorem in the free field
setting, see [\rfr(Kay)] for a short direct verification of this special case)
we may identify the rescaled boost transformations with the modular group of
the  real Hilbert subspace of the one-particle Hilbert space associated to the
corresponding wedge region. Then Property~(1.1) holds because it is equivalent
to the commutativity of  $\Theta$ with the boosts (see Proposition~2.6).
 \endProof
 \begremark{1.5.} The proof of Lemma~1.4 makes use of a one-particle version
of the Bisognano-Wichmann theorem; as we mentioned this can be proved directly
by mimicking the proof of the Bisognano-Wichmann
 in this special case. Since this quick verification requires an analytic
continuation argument that does not fit with the spirit of this paper, we
sketch
here an algebraic derivation of  Lemma~1.4. To begin with note that Lemma~1.4
would hold if  $V$ were the irreducible positive-energy massless representation
of $\Ppo(3)$ with helicity 0. Indeed in this case $V$ extends to a
representation
of the conformal group and the algebraic argument in [\rfr(BGL1)] applies. Now
$V\otimes V$ has a direct integral decomposition into irreducible
representations where massive representations occur and thus Lemma~1.2 implies
Lemma~1.4 for some $m>0$. However the representation $U\equiv V|_{\bG}$
in Lemma~1.4 does not depend on $m>0$ up to unitary equivalence by the
following
proposition, hence Lemma~1.4 holds for all $m>0$.
 \endremark
\begproposition{1.6.} The representation $U\equiv V|_{\bG}$
 in Lemma~1.4
is equivalent to the quasi-regular representation of $\bG$
corresponding to the rotation subgroup
 $\bK\equiv\left\{\left(\matrix{\cos\theta
&\sin\theta\cr-\sin\theta&\cos\theta\cr}\right),
0\leq\theta<2\pi\right\}/\{1,-1\}$,
 namely $U$ is the representation of $\bG$ induced by the identity
representation of $\bK$. The representation $U\otimes U$ is equivalent to an
infinite multiple of the regular representation $\lambda$ of $\bG$
 $$
U\otimes U = \infty\cdot\lambda.
 $$
 \endproposition
 \begProof The $m>0$ hyperboloid
 $H_m=\{\vec{x}\in\Re^3/x_0^2-x_1^2- x_2^2 = m^2,\, x_0>0\}$
 is a homogeneous space for $\bG$ whose stability subgroup at the point
$(m,0,0)$ is $\bK$. $U$ is the corresponding representation on
$L^2(H_m,\mu_m)$, with $\mu_m$ the Lorentz invariant measure on $H_m$, and this
is by definition the quasi-regular represention with respect to $\bK$.
 The last statement is a consequence of the Mackey tensor product
theorem for induced  representations, see [\rfr(Lipsman), Theorem~2 and
Example~5]. \endProof
 \begProofof{of Theorem~1.1} $(a)$ Property $(ii)$ for $\l$ follows by
Lemma~1.2
and Proposition~1.6. Here is an alternative verification of this fact. By
Proposition~1.6, taking tensor products and making use of Corollary~1.3, we
check that Property~(1.1) is valid for the irreducible representation $V_{m,s}$
of $\Ppo(3)$ of any mass $m>0$ and any integral spin $s$. Now  $U_{m,s}\equiv
V_{m,s}|_\bG$ is the induced representation
 $$
U_{m,s}=\Ind{\bK}{\bG}(\chi_s)
 $$
where $\chi_s$ is character of $\bK\simeq\To$ associated
with the integer $s$. By inducing at stages one has
 $$
\lambda=\Ind{\{1\}}{\bG}(id)=\Ind{\bK}{\bG}(\lambda_\bK)
 $$
where $\lambda_\bK=\Ind{\{1\}}{\bK}(id)$ is the
regular representation of \bK, hence
 $$
\lambda =\Ind{\bK}{\bG}(\lambda_\bK)
 = \Ind{\bK}{\bG}\left(\bigoplus_{s=-\infty}^{\infty}\chi_s\right) =
\bigoplus_{s=-\infty}^{\infty}\Ind{\bK}{\bG}(\chi_s) =
\bigoplus_{s=-\infty}^{\infty}U_{m,s}
 $$
 and the statement follows by Corollary~1.4.

\item{(b)} Let $U$ enjoy Property $(ii)$. If $U$ is a general representation of
$\widetilde\bG$, then it determines a projective representation of  $\bG$,
namely $U(l(\cdot))$ with $l$ a Borel section for the quotient map of
$\widetilde\bG$ modulo $\bG$. The tensor product $U\otimes\overline U$ of $U$
with  the conjugate representation is a true representation of $\bG$, hence
Property~$(ii)$ holds for $U\otimes\overline U$. Then, by Lemma~1.2,
Property~$(ii)$ holds for $U$.
 \endProof

\titlea{2. PCT, Spin and Statistics on the Field Algebras}

 In this section we consider a pre-cosheaf $\O\to\F(\O)$ of von~Neumann
algebras
acting on a Hilbert space $\H$, where $\O$ is any open region of
the $4$-dimensional Minkowski space $M$.
 We assume the following properties:
\medskip
\item{$(1)$} {\it Reeh-Schlieder property} for space-like cones:
there is a vector $\Q\in\H$ ({\it vacuum}) which is cyclic for the
algebras associated with all space-like cones.
\item{$(2)$} {\it Normal commutation relations}: there is a
vacuum-preserving self-adjoint unitary $\G$ (statistics operator)
that implements an automorphism on every local von Neumann algebra and the
normal commutation relations between Bose and Fermi fields hold, i.e. setting
 $$
\F_{\pm}(\O):=\{A\in\F(\O):\G A\G=\pm A\}
 $$
we have that if $\O_1$ and $\O_2$ are space-like separated then
$\F_+(\O_1)$ commutes with $\F(\O_2)$ and $\F_-(\O_1)$
anticommutes with $\F_-(\O_2)$.
\item{$(3)$} {\it Modular covariance property} with respect to the vacuum
vector $\Q$ (cf. Definition~2.3).
 \medskip
\begproposition{2.1}{\rm(Twisted Locality).} Let $Z$ be the unitary
operator defined by
 $$
Z={I+i\G\over 1+i}.\eqno(2.1)
 $$
 Then
 $$
Z\F(\O)Z^*\subset\F(\O')'.\eqno(2.2)
 $$
 \endproposition
 \begProof A direct computation shows that
 $$
\eqalign{
ZBZ^*=B,\qquad&B\in\F_+\cr
ZFZ^*=i\G F\quad&F\in\F_-\cr}.\eqno(2.3)
 $$
 Hence, if $\O_1$ and $\O_2$ are space-like separated and
$F_j\in\F_-(\O_j)$, $j=1,2$ we have
 $$
[ZF_1Z^*,F_2]=i\G(F_1F_2+F_2F_1)=0
 $$
and the thesis holds.
 \endProof
 We recall that a {\it wedge} region is any Poincar\'e
transformed of the region $W_1:=\{\vec{x}\in\Re^n:|x_0|<x_1\}$.
The {\it boosts} preserving $W_1$ are the elements of the
one-parameter subgroup $\L_{W_1}(t)$ of $\Ppo$ which acts
on the coordinates $x_0$, $x_1$ via the matrices
 $$
\left(\matrix{\cosh2\pi t&-\sinh2\pi t\cr
-\sinh2\pi t&\cosh2\pi t}\right)
 $$
 and leaves the other coordinates unchanged.
 The boosts $\L_W(t)$ for any wedge $W$ are defined by Poincar\'e
conjugation. We denote by $\W$ the family of all wedges in $M$.
 \par
 By twisted locality and the Reeh-Schlieder assumption the vacuum is cyclic and
separating for the algebras associated with all regions $\B$ such that both
$\B$
and $\B'$ contain some space-like cone. In particular, this property
holds for any wedge region $W$, therefore, by Tomita-Takesaki
theory, the modular operators $\D_W$ and $J_W$ associated with
$(\F(W),\Q)$ are defined for each region $W\in\W$.
 \beglemma{2.2.} The following commutation relations hold for all $W\in\W$:
 $$
\eqalign{
&[Z,\D_W^{it}]=0\cr
&J_WZJ_W=Z^*.\cr}
 $$
 \endlemma
 \begProof Since $\G$ preserves the vacuum and maps $\F(W)$ in
itself, it commutes with the modular operators. Then the relations
follow by the definition of $Z$.
 \endProof
 \begdefinition{2.3.} The pre-cosheaf $\O\to\F(\O)$, with the
Reeh-Schlieder property with respect to $\Omega$
 as above, satisfies {\it modular
covariance} if, for all regions $\O$,
 $$
\D^{it}_W\F(\O)\D^{-it}_W=\F(\L_W (t)\O),
\quad W\in\W\ .\eqno(2.4)
 $$
 \enddefinition
 \begproposition{2.4.} Let $\O\to\F(\O)$ be a modular covariant
pre-cosheaf of von~Neumann algebras satisfying the Reeh-Schlieder
property and normal commutation relations.
Then the modular unitaries associated to the wedges determine
a positive energy unitary representation $U$ of the covering
group $\Spo$ of the Poincar\'e group such that
 $$
U(g)\F(\O)U(g)^*=\F(\s(g)\O)\qquad g\in\Spo,\ \O\subset M
 $$
where $\s:\Spo\to\Ppo$ is the covering map, and
 $$
\D^{it}_W=U(\tilde\L_W(t)),\qquad W\in\W\eqno(2.5)
 $$
where $\tilde\L_W(t)$ denotes the lifting of $\L_W(t)$ in $\Spo$.
 \endproposition
 Before proving Proposition~2.4 we observe that twisted essential
duality holds if modular covariance is assumed.
 \begproposition{2.5}{\rm(Twisted essential duality).} If the pre-cosheaf
$\F$ satisfies modular covariance, then
 $$
Z\F(W)Z^*=\F(W')'.
 $$
As a consequence,
 $$
\D_W^{it}=\D_{W'}^{-it}.
 $$
 \endproposition
 \begProof By modular covariance, twisted locality and Lemma~2.2,
$Z\F(W)Z^*$ is a globally invariant subalgebra of $\F(W')'$  under
the action of the modular group $\D_{W'}^{it}$ for which the vacuum
is cyclic, therefore it coincides with $\F(W')'$. The last equality
follows by Lemma~2.2.
 \endProof
 \begProofof{of Proposition~2.4} The steps of the proof are essentially
the same as in the proof of Theorem~2.3~in~[\rfr(BGL2)], but
Lemma~2.5~[\rfr(BGL2)] should be reconsidered because we do not assume
additivity, and Fermi statistics is allowed. Indeed the first part, in which it
is shown that the set $H=\{g\in\Ppo:\F(g\O)=\F(\O)$\ ,\  $\O\subset M\}$ is a
normal subgroup of $\Ppo$, still holds.
 Therefore, with the same argument as in Lemma~2.6~[\rfr(BGL2)], we
get a central extension of $\Ppo/H$.

 When $H$ is $\{1\}$ or the translation subgroup, the triviality
of the central extensions for $\Ppo$ and for $\Lpo$ concludes the
theorem. When $H=\Ppo$ then the algebras associated to all wedges
coincide and, in particular, $\F(W)=\F(W')$. In this case
$\F_+(W)$ is abelian and two Fermi fields localized in $W$
anticommute. Then if $A\in\F_-(W)$ we have $A^*A+AA^*=0$ which
implies $\F_-(W)=0$. As a consequence $\F(W)$ is abelian for
any wedge $W$ and the representation of $\Spo$ implemented by the
modular unitaries is the trivial representation.
 \endProof
 We notice that in the proof of the preceding theorem the Reeh-Schlieder
property was needed only for wedge regions. When this property does not hold
for
double cones it is possible to construct theories for which the group $H$ in
the
proof is the translation group, namely the translations act identically on the
pre-cosheaf (see [\rfr(GL3)]).

  Our next step is to show that, as a consequence of Theorem~1.1,
the representation $U$ extends to (some covering of) the proper
Poincar\'e group $\Pp$.
 \begproposition{2.6.} Let $W_1$, $W_2$ be orthogonal wedges,
and let us denote by $\D_i$, $J_i$, $i=1,2$ the modular operators
of $(\F(W_i),\Q)$. Then the following relation holds
 $$
\D_1^{1/2}\D_2^{it}\D_1^{-1/2}\subset J_1\D_2^{it}J_1\ ,\eqno(2.6)
 $$
 and the left hand side is densely defined.
 \endproposition
 \begProof We may suppose that $W_i:=\{\vec{x}\in\Re^n:|x_0|<x_i\}$,
$i=1,2$.
 With the same arguments of Lemma~1.2, it is enough to show
that, for each $t\in\Re$, there exists a dense set $\cD_t$ such
that
 $$
\D_1^{1/2}\D_2^{it}\D_1^{-1/2}\xi=J_1\D_2^{it}J_1\xi
\qquad\xi\in\cD_t\ .
 $$
 We claim that we may take
 $$
\cD_t=\{J_1A\Q:A\in\F(W_1\cap\L_2(-t)W_1)\}\ .
 $$
 First we observe that Relation~(2.6) holds on $\cD_t$. Indeed, for
$\xi=J_1A\Q\in\cD_t$,
 $$
\eqalign{
J_1\D_2^{it}J_1\xi
&=J_1\D_2^{it}A\Q=J_1(\D_2^{it}A\D_2^{-it})\Q
 =\D_1^{1/2}\D_2^{it}A^*\Q\cr
&=\D_1^{1/2}\D_2^{it}S_1J_1\xi
 =\D_1^{1/2}\D_2^{it}\D_1^{-1/2}\xi\ .\cr}
 $$
 Now, by the Reeh-Schlieder property for space-like cones, it is
enough to show that $\L_2(-t)W_1\cap W_1$ contains a space-like
cone in order to prove that $\cD_t$ is dense. Indeed the
space-like cone
 $$
\C_t=\left\{\vec{x}\in M:
\sqrt{x_0^2+x_2^2}<{x_1\over\sqrt2\cosh2\pi t}\right\}
 $$
is obviously contained in $W_1$ and, if
$\vec{y}=\L_2(t)\vec{x}$, $\vec{x}\in\C_t$, then
 $$
|y_0|=|\cosh2\pi t x_0-\sinh2\pi t x_2|
\leq\sqrt2\cosh2\pi t\sqrt{x_0^2+x_2^2}<x_1\ ,
 $$
namely $\L_2(t)\C_t\subset W_1$.
 \endProof
 \begcorollary{2.7.} Let $W_1$, $W_2$ be orthogonal wedges, and
let us denote by $\D_i$, $J_i$, $i=1,2$ the modular operators of
$(\F(W_i),\Q)$. Then, for each $t\in\Re$, the following relation
holds:
 $$
J_1\D_2^{it}J_1=\D_2^{-it}\ .\eqno(2.7)
 $$
 \endcorollary
 \begProof By Proposition~2.6 it is enough to show that, for each
$t\in\Re$,
 $$
\D_1^{1/2}\D_2^{it}\D_1^{-1/2}\subset\D_2^{-it}\ .
 $$
 Since the two wedges are orthogonal, the homomorphism
$\f:{\rm SL}(2,\Re)\to\Spo$ determined by $\f(\m(t))=\tilde\L_{W_1}(t)$ and
$\f(\n(t))=\tilde\L_{W_2}(t)$ is injective.
Then we apply Theorem~1.1 and the proof is complete.
 \endProof
 Let us consider the automorphism $\a$ of the Poincar\'e group
given by $\a(g)=I_1gI_1$, where $I_1$ is the diagonal matrix which
changes the sign of the $x_0$ and $x_1$ components and leave the
other coordinates fixed. The proper Poincar\'e group $\Pp$
is generated by $I_1$ and $\Ppo$, as a consequence it may be seen
as a semidirect product:
 $$
\Pp=\Ppo\times_\a\Ze_2.\eqno(2.8)
 $$
 The action $\a$ gives rise to an action on the Lie algebra of
$\Ppo$ and therefore to an action on the universal covering
$\Spo$ of $\Ppo$.
 \par
 Since the center of $\Spo$ has order $2$, $\a$
is trivial on the center.
 \par
 Now let us denote by $\Sp$ the semidirect
product of the universal covering of $\Ppo$ with $\Ze_2$ via $\a$.
It follows that $\Sp$ is a central extension of the proper
Poincar\'e group $\Pp$ via the homomorphism $\s$ extended by
$\s(\tilde I_1)=I_1$, where $\tilde I_1$ denotes the non-trivial
element in the $\Ze_2$ component of $\Sp$.
 \begproposition{2.8.} The representation $U$ extends to a
(anti-)unitary representation (still denoted by $U$) of the group
$\Sp$ with $U(\tilde I_1)=ZJ_1$.
 \endproposition
 \begProof We only have to show that the following relation holds:
 $$
ZJ_1U(g)(ZJ_1)^*=U(\a(g)),\qquad g\in\Spo.\eqno(2.9)
 $$
Relation~(2.9) holds for $\tilde\L_{W_1}(t)$ by Tomita-Takesaki
theory, for the translations by Borchers theorem [\rfr(Borc1)] and
for $\tilde\L_{W}(t)$ when $W$ is a wedge orthogonal to $W_1$ by
Lemma~2.2 and Proposition~2.5.
Since the lifting of the boosts associated with a maximal set of
orthogonal wedges and the translations generate $\Spo$, we get
the the\-sis.
 \endProof
 We may associate a dual pre-cosheaf $\F^d$ to the given pre-cosheaf
$\F$. If $\O$ is a convex causally complete region the algebra
$\F^d(\O)$ is defined by
 $$
\F^d(\O)\equiv\bigcap_{W\supset\O}\F(W),
 $$
 where $W$ is a wedge. For any more general region $\B$ we set
 $$
\F^d(\B)\equiv\left(\bigcup_{\O\subset\B}\F(\O)\right)'',
 $$
where $\O$ is a convex causally complete region. It is easy to see
that $\F^d$ is a pre-cosheaf satisfying Reeh-Schlieder property,
locality and duality for convex causally complete regions.

 \begproposition{2.9.} The representation $U$ of $\Sp$
implements a covariant action on the wedges, i.e.
 $$
U(g)\F(W)U(g)^*=\F(gW)\ ,\qquad\forall g\in\Sp,\ W\in\W\ .
 $$
 As a consequence, the action is covariant on the dual pre-cosheaf
$\O\to\F^d(\O)$.
 \endproposition
 \begProof
By Proposition~2.4, we have to prove covariance only for the element
$ZJ_1$.
 For each $W\in\W$ we may find a $g\in\Sp$ such that $W=\s(g)W_1$.
 Then
 $$
\eqalign{
ZJ_1\F(W)(ZJ_1)^*&=ZJ_1U(g)\F(W_1)(ZJ_1U(g))^*\cr
&=ZJ_1U(g)ZJ_1\F(I_1W_1)(ZJ_1U(g)ZJ_1)^*\cr
&=U(\a(g))\F(-W_1)U(\a(g))^*\cr
&=\F(I_1W)\ .\cr}
 $$
 Since the dual pre-cosheaf $\F^d$ is described in terms of the
restriction of $\F$ to the wedges, the second part of the statement
follows.
 \endProof
 \begtheorem{2.10}{\rm(PCT).} Let $\O\to\F(\O)$ be a modular covariant
pre-cosheaf satisfying Reeh-Schlieder property and normal commutation relations
as above. Then there exists an antiunitary operator $\T$ which implements the
PCT symmetry on the dual pre-cosheaf, i.e.
 $$
\T\F^d(\O)\T=\F^d(-\O)\ .\eqno(2.10)
 $$
 \endtheorem
 \begProof Since the dimension of the space-time is even, the transformation
which changes the sign of all coordinates belongs to the proper Poincar\'e
group
$\Pp$. Then we choose a preimage $\theta\in\Sp$ via $\s$ of this element
and set $\T:=U(\theta)$. Obviously $\T$ verifies Eq.~(2.10)
and is anti-unitary.
 \endProof
 We remark that Proposition~2.9 is an abstract form of the PCT theorem
and is indeed what we can prove in an odd dimensional Minkowski space, where
the global space-time inversion does not belong to the proper Poincar\'e group.

 \begtheorem{2.11}{\rm(Spin and Statistics).} Let $\O\to\F(\O)$
be a modular covariant pre-cosheaf of field algebras as in Theorem~2.10. Then
 $$
\G=U(2\pi).\eqno(2.11)
 $$
 \endtheorem
 \begProof Let $W_2$ be the wedge along the axis $x_2$. It is well known that
the generator of the rotations $r(\th)$ on the $x_1-x_2$ plane and such that
$r(\pi/2)W_1=W_2$ is a multiple of the commutator between the generators of
$\L_1(t)$ and $\L_2(t)$. Then, by definition, $\a(\tilde r(\th))=\tilde r(-
\th)$,
where $\tilde r$ denotes the lifting of $r$ in $\Spo$. Therefore we have
 $$
ZJ_1ZJ_2=ZJ_1U(\tilde r({\pi\over2}))ZJ_1U(\tilde r(-{\pi\over2}))=
U(\a(\tilde r({\pi\over2})))U(\tilde r(-{\pi\over2}))=U(\tilde
r(-\pi))\ .
 $$
 On the other hand, since $I_1W_2=W_2$, Proposition~2.9 shows that
$ZJ_1$ leaves $\F(W_2)$ globally stable, hence it commutes
with $J_2$. Therefore
 $$
ZJ_1ZJ_2=ZJ_1J_2Z^*=J_2ZJ_1Z^*=Z^*J_2ZJ_1=\G ZJ_2ZJ_1\ .
 $$
 Finally,
 $$
U(2\pi)=(ZJ_1ZJ_2)^2=(\G ZJ_2ZJ_1)(ZJ_1ZJ_2)=\G\ .
 $$
 \endProof

\titlea{3. PCT and Spin and Statistics Theorems for
Superselection Sectors with Finite Statistics}

 In this section we shall prove the Spin and Statistics theorem in
the framework of local observable algebras [\rfr(Haag1)], i.e. we
shall prove that a  para-Bose (resp. para-Fermi) superselection sector with
finite statistics is covariant with respect to the group  $\Ppo$ (resp.
$\Spo$).

 For the sake of simplicity, in this section we shall make the usual assumption
the von Neumann algebras associated with unbounded regions are generated
by additivity by the ones associated with double cones
within the region,   so that the Doplicher-Roberts theorem on the
reconstruction
of the field algebras applies [\rfr(DoRo)].

We shall consider morphisms of the local observable algebras localized in
space-like cones in a 4-dimensional space-time [\rfr(DHR1),\rfr(BuFr1)]. The
same techniques  would prove the Spin and Statistics theorem for sectors
localized in space-like cones  in the $n$-dimensional Min\-kow\-ski space for
any $n\geq4$ and for  sectors localized in double-cones if $n=3$. We omit here
the treatment of these cases in view of a  more general intrinsic proof
[\rfr(GL3)], valid also in the low-dimensional case, that will be carried on in
the setting of the
 of pre-cosheaves of von Neumann algebras on wedge regions.

 In the following we  consider a net of local observable von Neumann algebras
[\rfr(HaKa1)], i.e.  a map
 $$
\A:\O\to\A(\O),\qquad \O\in\K
 $$
where $\K$ is the family of the double cones in the
$4$-dimensional Minkowski space $M$, such that
 $$
\eqalignno{ \O_1\subset\O_2&\imply\A(\O_1)\subset\A(\O_2)
& {\rm(isotony)}\cr
\A(\O)&\subset\A(\O')'& {\rm(locality)}.\cr}
 $$
 \par
 If $\B$ is an unbounded open region of $M$, we shall denote by
$\A(\B)$ the von Neumann algebra generated by the
von~Neumann algebras $\A(\O)$, $\O\subset\B$, so that $\A$ extends to a
pre-cosheaf of von Neumann algebras on more general regions. The local algebras
are supposed to act on a separable Hilbert space $\H_0$, with a common cyclic
vector $\Q_0$.
  \par
The modular covariance property with respect to $\Q_0$ is  assumed.
 As a particular case of the analysis in the previous section, this property
gives rise to a positive-energy representation of the Poincar\'e group which
leaves $\Q_0$ invariant. By locality $\Q_0$
is cyclic and separating for the von Neumann algebras associated with all
non-empty open regions whose complement have non-empty interior (Reeh-Schlieder
property).
 \begproposition{3.1.} $\A$ has a unique direct integral decomposition into
irreducible modular covariant nets.
 \endproposition
 \begProof
 If $\xi\in\H_0$ is a $\D_W^{it}$--invariant vector for some wedge $W$,
then $\xi$ is also $\D_{W_1}^{it}$--invariant for all other wedges $W_1$. This
follows by the representation theory of ${\rm SL}(2,\Re)$, since the modular
operators generate a representation of the Poincar\'e group. Therefore the
centralizer $\Z$ of $\A(W)$ with respect to $\Omega_0$ does not depend on $W$.
By locality $\Z$ is an abelian von Neumann algebra. The direct integral
decomposition with respect to $\Z$ has then the desired properties. The
uniqueness is clear since $\Z$ is canonically constructed.
 \endProof
Because of Proposition~3.1 we shall always assume in the following that the net
is irreducible.
 We summarize now the results so far obtained in the following theorem.
 \begtheorem{3.2.} Let $\O\to\A(\O)$ an irreducible local, modular
covariant net of von~Neumann algebras on $M$ as above. The extended
pre-cosheaf satisfies essential duality, and there is a (unique) positive
energy unitary representation $g\in\Ppo\to U_0(g)$ with the
following properties:
 \medskip
 \item{$(a)$} $U_0(g)\A(\O)U_0(g)^*=\A(g\O)$,
 \item{$(b)$} $U_0(\L_W(t))=\D_W^{it}$, $W\in\W$,
 \item{$(c)$} {\rm[PCT Theorem]} There is an antiunitary operator
$\T$ (PCT operator) satisfying
 $$
\eqalign{
\T\A(\O)\T&=\A(-\O),\qquad\O\in\K\cr
[{\rm ad}\T\cdot\r]&=[\overline\r],\cr}
 $$
where $\r$ is a morphism localized in a space-like cone, and  the PCT operator
is given by $\T=J_WR_W$ where $W$ is any wedge whose edge contains the origin
and $R_W$ is the space reflection w.r.t. all directions contained in the edge
of
$W$.
 \medskip
 \endtheorem
 \begProof Since the net $\A$ may be considered as a purely Bose
pre-cosheaf of field algebras (cf. Sect.~2), the thesis
follows by Theorems~2.9 and~2.10 and by [\rfr(GL1)].
 \endProof
 The previous theorem shows in particular that the Spin and Statistics theorem
holds in the vacuum sector,  solving a problem that remained open in
[\rfr(BGL2)].
 \bigskip
 We recall now that a localized representation of the quasi-local
C$^*$-algebra
 $$
\A := \left(\bigcup_{\cal O\in K}\A(\cal O)\right)^-
 $$
 (norm closure) is a representation which is equivalent to the vacuum
representation when restricted to the space-like complement of any space-like
cone.
 Such a representation is equivalent a localized
transportable morphism, that is to say a morphism $\r$ of $\A$
into $B(\H_0)$ which is localized in a given space-like cone $\cal S$,
i.e.
$\r(x)=x$ for each $x\in\A(\cal S')$.

 A unitary equivalence class of localized representations is
called a {\it superselection sector}.

 A localized morphism $\r$ is called Poincar\'e covariant if
there exists a unitary representation $U_\r$ of the covering of the
Poincar\'e group $\Spo$ verifying
 $$
U_\r(g)\r(x)U_\r(g)^*=\r(U_0(\s(g))xU_0^*(\s(g))),
\qquad g\in\Spo,\ x\in\A.
 $$

 We recall that when the net satisfies a mild condition called
{\it regularity}, all localized transportable morphisms  with
finite statistics are Poincar\'e covariant  [\rfr(GL1)].

 As already mentioned, a sector carries a projective representation
of the Poincar\'e group, i.e. $U_\r(2\pi)$ is not necessarily
the identity. If $\r$ is irreducible $U_\r(2\pi)=\pm1$.

If an irreducible morphism $\r$ with finite statistics is
localized in a double cone, Doplicher, Haag and Roberts [\rfr(DHR1)] showed
that  the statistical behaviour of $\r$ is described by the
statistical parameter $\l_\r$, an analysis extended to topological charges by
Buchholz and Fredenhagen [\rfr(BuFr1)]. As explained in the Introduction, the
spin-statistics relation takes the form of an equality of complex numbers, i.e.
$U_\r(2\pi)=\rm{sign}(\l_\r)$.

 \begtheorem{3.3}{\rm(Spin and Statistics).} Let $\O\to\A(\O)$ be an
irreducible, local, modular covariant net of von~Neumann algebras on the
Minkowski space and let $\r$ be an irreducible, covariant, localized morphism
with finite statistics.
 Then
 $$
U_\r(2\pi)=\rm{sign}(\l_\r).
 $$
\endtheorem

By the Doplicher-Roberts theorem [\rfr(DoRo)], the superselection structure is
described by a field net $\O\to\F(\O)$ of von~Neumann algebras satisfying
normal
commutation relations and by a vacuum-preserving unitary representation $V$ of
a
compact group $G$ (the gauge group) on a Hilbert space $\H$, that implements
automorphisms on any $\F(\O)$.

 There is a canonical net isomorphism $\p:\A(\O)\to\F^G(\O)$ where $\F^G(\O)$
is the fixed-point algebra under the gauge action, and the covariance
property of the sectors gives a uniquely determined representation $U$ of
$\Spo$
on $\H$ which commutes with the gauge group representation and  restricts to
$U_0$ on the vacuum Hilbert space $\H_0$.

 \beglemma{3.4.} The representation $U$ satisfies
 $$
U(\L_W(t))=\D_W^{it},\qquad W\in\W,\ t\in\Re\ .
 $$
 \endlemma
 \begProof With $W$ a wedge region, denote by $\Delta_W$ the modular operator
associated with $\F(W)$.
Observe that $c_W(t):=U(\L_W(t))\D_W^{-it}$ implements an
automorphism of the algebra $\F(W)$  leaving
$\F^G(W)$ pointwise fixed by Theorem~3.1$(b)$.
The translations along $W$ commute with $c_W(t)$, since they commute both
with $U(\L_W(t))$and $\D_W^{it}$, and act ergodically on $\F(W)$ by the
irreducibility of the net.
By the Galois correspondence [\rfr(AHKT1)], it follows that $c_W(\,\cdot\,)$ is
a one-parameter subgroup of $V(G)$.

 We show now that $c_W(t)=I$.
 Indeed, on the one hand, by twisted duality (Proposition~2.5)
 $$
c_{W'}(t)=c_W(-t),\qquad W\in\W,\ t\in\Re.
 $$

 On the other hand $c_{gW}(t)=U(g)c_{W}(t)U(g)^*=c_{W}(t)$
because the gauge group action commutes with the representation $U$ of
$\Spo$, which means that $c_{W}(t)$ does not depend by $W$. Then
$c_W(t)=I$, i.e. $U(\L_W(t))=\D_W^{it}$.
\endProof

 \begProofof{of Theorem~3.3}
Each irreducible superselection sector   $\rho$ is associated in a one-to-one
correspondence with an irreducible representation $V_\rho$ of the gauge
group $G$. The Hilbert space $\H$  decomposes into a direct sum
 $$
\H=\bigoplus_{\r}\H_\r,
 $$
 where $\H_\r:=\H_0\otimes\Co^{d(\r)}$ and $d(\r)=|\l_\r|^{-1}= {\rm
dim}(V_\r)$ is the statistical dimension of $\r$, so that the
representation $\pi$ of the observable algebra $\cal A$ on $\H$, the
representation $U$ of $\Spo$ and the representation $V$ of
the gauge group decompose  accordingly as
 $$
\eqalign{
\pi(x)&=\bigoplus_{\r\in\hat G}\r(x)\otimes I,\quad x\in\A,\cr
U(g)&=\bigoplus_{\r\in\hat G}U_\r(g)\otimes I,\quad g\in\Spo,\cr
V(g)&=\bigoplus_{\r\in\hat G}I\otimes V_\r(g),\quad g\in G,\cr}
 $$
 where $U_\r$ is the representation of $\Spo$ in the sector $\r$
as before.

 In particular the statistics operator $\G$ is constant on the
direct summands $\H_\r$, and we have
 $$
\G=\bigoplus_{\r}\k_\r I_{\H_\r}
 $$
where $\k_\r=\rm{sign}(\l_\r)$ is the statistics phase.

 Then the Spin and Statistics theorem  follows directly by
Theorem~2.10.
 \endProof

\acknow{We thank Klaus Fredenhagen for pointing out a gap in the proof of
Theorem~1.1 of our original manuscript.}

\references

\end